\documentclass[12pt]{article}
\usepackage{graphicx,amssymb}
\usepackage{graphics}
\usepackage{epsfig}
\usepackage{subfigure}
\usepackage{authblk}
\usepackage{fullpage}
\usepackage{amsmath}

\begin{document}

\title{Combined prefactored compact schemes for first- and second-order derivatives: conceptual derivation}

\author[]{\normalsize  Adrian Sescu \thanks{sescu@ae.msstate.edu} }

\affil[]{Department of Aerospace Engineering, Mississippi State University}

\date{}

\maketitle

\begin{abstract}

The derivation of combined prefactored compact schemes for first and second order derivatives is described here, relying on the Fourier analysis of the original prefactored compact schemes. By this approach, the order of accuracy of the original schemes can be increased from sixth to eight, or from eight to tenth (depending on the order of the original scheme), while the number of grid points in the stencil is kept the same. Here, we only frame the conceptual derivation of the schemes, leading to a closed set of equations for the weights.

\end{abstract}

\maketitle

\section{Introduction}

Compact difference schemes possess the advantage of attaining higher-order of accuracy with fewer grid points per stencil. They are preferred in applications where high accurate results are desired, such as direct numerical simulations, large eddy simulations, computational aeroacoustics or electromagnetism, to enumerate few. One of the disadvantages of compact schemes is that an implicit approach is required to determine the grid functions, wherein a large matrix has to be inverted. A comprehensive study of high-order compact schemes approximating both first and second derivatives on a uniform grid was performed by Lele \cite{lele}. A wavenumber based optimization was introduced wherein the dispersion error was reduced significantly, achieving spectral-like resolution. Over the next years, compact schemes have been studied by many research groups, and applied to various engineering problems (see for example, Li et al. \cite{Li}, Adams and Shariff \cite{Adams}, Liu \cite{Liu}, Deng and Maekawa \cite{Deng1}, Fu and Ma \cite{Fu1,Fu2}, Meitz and Fasel \cite{Meitz}, Shen et al. \cite{Shen}, Shah et al. \cite{Shah}). Other worth-mentioning examples include Kim and Lee \cite{Kim} who performed an analytic optimization of compact finite difference schemes, Mahesh \cite{mahesh} who derived a family of compact finite difference schemes for the spatial derivatives in the Navier-Stokes equations based on Hermite interpolations (see also, Chu and Fan \cite{chu} for a similar prior analysis), or Deng and Zhang \cite{Deng2} who developed compact high-order nonlinear schemes which are equivalent to fifth-order upwind biased explicit schemes in smooth regions. 

Hixon \cite{hixon1,hixon2} derived prefactored high-order compact schemes that use three-point stencils and returns up to eighth-order of accuracy. These schemes combine the tridiagonal compact formulation with the optimized split derivative operators of an explicit MacCormack type scheme. The optimization of Hixon's \cite{hixon1,hixon2} schemes in terms of reducing the dispersion error was performed by Ashcroft and Zhang \cite{Ashcroft} who used Fourier analysis to select the coefficients of the biased operators such that the dispersion characteristics match those of the original centered compact scheme and their numerical wavenumbers have equal and opposite imaginary components. Today, compact schemes are widely used in numerical simulations of turbulent flows (e.g., direct numerical simulations), computational aeroacoustics, or computational electromagnetics. In order to increase the speed of such numerical simulations it is desirable to derive more computational efficient compact schemes without affective the order of accuracy and the wavenumber characteristics.

In this work, we frame the conceptual derivation of combined prefactored compact schemes for first and second order derivatives, aimed at increasing the resolution accuracy with fewer points per stencil. They are based on the type of prefactorization introduced previously by Hixon \cite{hixon1,hixon2}. One of the advantages of these schemes is that the derivatives are explicitly determined by sweeping from one boundary to the other, thus avoiding the inversion of matrices which can increase the computational time significantly.

\section{Derivation of the combined prefactored schemes}

Chu and Fan~\cite{chu} derived a combined compact difference scheme (over a three-point stencil), which is sixth-order accurate (see also Mahesh~\cite{mahesh} for a generalization). The scheme consists of two coupled equations for the first and second derivatives and can be written as:

\begin{eqnarray}\label{1}
u_i'+ \frac{7}{16}
\left(
u_{i+1}' + u_{i-1}'
\right) -
\frac{h}{16}
\left(
u_{i+1}'' - u_{i-1}''
\right) =
\frac{15}{16h}
\left(
u_{i+1} - u_{i-1}
\right)
\end{eqnarray}

\begin{eqnarray}\label{2}
u_i''+ \frac{9}{8h}
\left(
u_{i+1}' - u_{i-1}'
\right) -
\frac{1}{8}
\left(
u_{i+1}'' + u_{i-1}''
\right) =
\frac{3}{h^2}
\left(
u_{i+1} + u_{i-1}
\right) -
\frac{6}{h^2}u_{i}
\end{eqnarray}
or, following the notation used in~\cite{hixon1},

\begin{eqnarray}\label{3}
D_i+ \frac{7}{16}
\left(
D_{i+1} + D_{i-1}
\right) -
\frac{h}{16}
\left(
D_{i+1}^2 - D_{i-1}^2
\right) =
\frac{15}{16h}
\left(
u_{i+1} - u_{i-1}
\right)
\end{eqnarray}

\begin{eqnarray}\label{4}
D_i^2+ \frac{9}{8h}
\left(
D_{i+1} - D_{i-1}
\right) -
\frac{1}{8}
\left(
D_{i+1}^2 + D_{i-1}^2
\right) =
\frac{3}{h^2}
\left(
u_{i+1} + u_{i-1}
\right) -
\frac{6}{h^2}u_{i}
\end{eqnarray}
where $h$ is the grid step, $u$ is the grid function, and $D$ and $D^2$ stand for first and second order derivatives. Mahesh~\cite{mahesh} extended the idea to higher-order, such as the eight-order coupled stencils for the first and second derivatives in the form:

\begin{eqnarray}\label{5}
D_i+ \frac{17}{36}
\left(
D_{i+1} + D_{i-1}
\right) -
\frac{h}{12}
\left(
D_{i+1}^2 - D_{i-1}^2
\right) =
\frac{107}{108h}
\left(
u_{i+1} - u_{i-1}
\right) -
\frac{1}{108h}
\left(
u_{i+2} - u_{i-2}
\right)
\end{eqnarray}

\begin{eqnarray}\label{6}
D_i^2 +
\frac{23}{18h}
\left(
D_{i+1} - D_{i-1}
\right) -
\frac{1}{6}
\left(
D_{i+1}^2 + D_{i-1}^2
\right) =
\frac{88}{27h^2}
\left(
u_{i+1} + u_{i-1}
\right) -
\frac{1}{108h^2}
\left(
u_{i+2} - u_{i-2}
\right) -
\frac{13}{2h^2}
u_{i}
\end{eqnarray}

Equations (\ref{3}) and (\ref{4}) or (\ref{5}) and (\ref{6}) can be jointly written in matrix form:

\begin{eqnarray}\label{7}
[B]\{D\} = [C]\{u\}
\end{eqnarray}
where $[B]$ and $[C]$ are $2N\times 2N$ matrices, $\{D\}$ is a vector of $2N$ elements containing the first and the second derivatives for $i=1,2,...,N$ and $\{u\}$ is a vector of $2N$ elements containing the dependent variables twice ($N$ is the number of grid points along a grid line).

To derive the prefactored compact schemes, the forward and backward operators $D_i^F$ and $D_i^B$ for the first derivative and $D_i^{2F}$ and $D_i^{2B}$ for the second derivative are defined as

\begin{eqnarray}\label{8}
D_i = \frac{1}{2}(D_i^F + D_i^B)
\end{eqnarray}
and

\begin{eqnarray}\label{9}
D_i^2 = \frac{1}{2}(D_i^{2F} + D_i^{2B})
\end{eqnarray}

The spatial stencils for the forward and backward derivative operators are defined as a combination of both the first and second derivatives:

\begin{eqnarray}\label{10}
D_i^F +
\beta_F^I D_{i+1}^F +
\theta_F^I h D_{i+1}^{2F} =
\frac{1}{h}
\left(
a_F^I u_{i-1} + b_F^I u_{i} + c_F^I u_{i+1}
\right)
\end{eqnarray}

\begin{eqnarray}\label{11}
D_i^B +
\beta_B^I D_{i-1}^B +
\theta_B^I h D_{i-1}^{2B} =
\frac{1}{h}
\left(
a_B^I u_{i-1} + b_B^I u_{i} + c_B^I u_{i+1}
\right)
\end{eqnarray}
and

\begin{eqnarray}\label{12}
D_i^{2F} +
\beta_F^{II} \frac{1}{h} D_{i+1}^F +
\theta_F^{II} D_{i+1}^{2F} =
\frac{1}{h^2}
\left(
a_F^{II} u_{i-1} + b_F^{II} u_{i} + c_F^{II} u_{i+1}
\right)
\end{eqnarray}

\begin{eqnarray}\label{13}
D_i^{2B} +
\beta_B^{II} \frac{1}{h} D_{i-1}^B +
\theta_B^{II} D_{i-1}^{2B} =
\frac{1}{h^2}
\left(
a_B^{II} u_{i-1} + b_B^{II} u_{i} + c_B^{II} u_{i+1}
\right)
\end{eqnarray}
where the coefficients $a$, $b$, and $c$ must be chosen such that when the two biased stencils are added, the original stencils of Chu and Fan~\cite{chu} or Mahesh~\cite{mahesh} are recovered. Following the analysis of Hixon and Turkel~\cite{hixon1}, the real components of the numerical wavenumbers of the forward and backward stencils are required to be equal and identical to the numerical wavenumber of the original scheme, and the imaginary components of the numerical wavenumbers are required to be equal and opposite.

To apply the Fourier analysis, the dependent variable $u$ is assumed to be periodic in the domain $[0,L]$. The Fourier decomposition of $u$ is:

\begin{eqnarray}\label{14}
u(x) =
\sum_{k=-N/2}^{k=N/2}
\hat{u}_k e^{\frac{2\pi ikx}{L}}
\end{eqnarray}
where $i=\sqrt{-1}$. To simplify the analysis, a scaled wavenumber $w=2\pi kh/L$ and a scaled coordinate $s=x/h$ are introduced, such that the Fourier modes are simply $exp(iws)$. The first and the second derivatives of the exact Fourier coefficients are given by:

\begin{eqnarray}\label{15}
\hat{u}_k' =
\frac{iw}{h}\hat{u}_k,
\hspace{6mm}
\hat{u}_k'' =
-\left(
\frac{w}{h}
\right)^2
\hat{u}_k
\end{eqnarray}
while the Fourier coefficients of the derivatives obtained from the differencing schemes are:

\begin{eqnarray}\label{16}
(\hat{u}_k')_{CCD} =
\frac{iw'}{h}\hat{u}_k,
\hspace{6mm}
(\hat{u}_k'')_{CCD} =
-\left(
\frac{w''}{h}
\right)^2
\hat{u}_k
\end{eqnarray}
where $w'=w'(w)$ and $w''=w''(w)$ are the modified wavenumbers for the first and second order derivatives, respectively.

The numerical wavenumbers of the original combined compact scheme are given by~\cite{chu}

\begin{eqnarray}\label{17}
w'(w) = 
\frac{9\sin{w}[4+\cos{w}]}{24+20\cos{w}+\cos{2w}}
\end{eqnarray}
and

\begin{eqnarray}\label{18}
w''^2(w) = 
\frac{81-48\cos{w}-33\cos{2w}}{48+40\cos{w}+2\cos{2w}}
\end{eqnarray}
for stencils (\ref{3}) and (\ref{4}) (Chu and Fan~\cite{chu}) and

\begin{eqnarray}\label{19}
w'(w) = 
\frac{\sin{w}[293+126\cos{w}+\cos{2w}]}
{6(34+33\cos{w}+3\cos{2w})}
\end{eqnarray}
and

\begin{eqnarray}\label{20}
w''^2(w) = 
\frac{1730-675\cos{w}-10870\cos{2w}-29\cos{3w}}
{36(34+33\cos{w}+3\cos{2w})}
\end{eqnarray}
for stencils (\ref{5}) and (\ref{6}) (Mahesh~\cite{mahesh}).

For the prefactored CCD schemes the modified wavenumbers $w_F'(w)$, $w_B'(w)$, $w_F''(w)$ and $w_B''(w)$ can be determined from the equations obtained by applying the Fourier transform to (\ref{10})-(\ref{13}):

\begin{eqnarray}\label{21}
iw_F' +
i\beta_F^I  w_F' e^{iw} -
\theta_F^I  (w_F'')^2 e^{iw} =
a_F^I e^{-iw} + b_F^I + c_F^I e^{iw}
\end{eqnarray}

\begin{eqnarray}\label{22}
- (w_F'')^2 +
i\beta_F^{II}  w_F' e^{iw} -
\theta_F^{II}  (w_F'')^2 e^{iw} =
a_F^{II} e^{-iw} + b_F^{II} + c_F^{II} e^{iw}
\end{eqnarray}
for the forward operators and

\begin{eqnarray}\label{23}
i w_B' +
i\beta_B^I  w_B' e^{-iw} -
\theta_B^I  (w_B'')^2 e^{-iw} =
a_B^I e^{-iw} + b_B^I + c_B^I e^{iw}
\end{eqnarray}

\begin{eqnarray}\label{24}
- (w_B'')^2 +
i\beta_B^{II}  w_B' e^{-iw} -
\theta_B^{II}  (w_B'')^2 e^{-iw} =
a_B^{II} e^{-iw} + b_B^{II} + c_B^{II} e^{iw}
\end{eqnarray}
for the backward operators. Equations (\ref{21}) and (\ref{22}) can be written in matriceal form as

\begin{eqnarray}\label{25}
[A]\{X\} = \{R\}
\end{eqnarray}
where $[A]$ is a $4\times 4$ matrix given by

\begin{eqnarray}\label{26}
[A] = \left[ \begin{array}{ccccc}
-\beta_F^I\sin{w} & -1-\beta_F^I\cos{w} & -\theta_F^I\cos{w}     & \theta_F^I\sin{w}  \\
1+\beta_F^I\cos{w} & -\beta_F^I\sin{w} & -\theta_F^I\sin{w}     & -\theta_F^I\cos{w}  \\
-\beta_F^{II}\sin{w} & -\beta_F^{II}\cos{w} & -1-\theta_F^{II}\cos{w}     & \theta_F^{II}\sin{w}  \\
\beta_F^{II}\cos{w} & -\beta_F^{II}\sin{w} & -\theta_F^{II}\sin{w}     & -1-\theta_F^{II}\cos{w}  \end{array} \right] ,
\end{eqnarray}
and $\{X\}$ is the unknown vector of 4 elements given by

\begin{eqnarray}\label{27}
\{X\} =
\{ \begin{array}{c}
\Re{(w_F')}
\hspace{4 mm}
\Im{(w_F')}
\hspace{4 mm}
\Re{((w_F'')^2)}
\hspace{4 mm}
\Im{((w_F'')^2)}
  \end{array} \}^T,
\end{eqnarray}
where $\Re()$ and $\Im()$ stands for real and imaginary parts, respectively. $\{R\}$ in (\ref{25}) is the right hand side vector given by

\begin{eqnarray}\label{28}
\{R\} =
\{ \begin{array}{c}
(c_F^I+a_F^I)\cos{w}+b_F^I
\hspace{4 mm}
(c_F^I-a_F^I)\sin{w}
\hspace{4 mm}
(c_F^{II}+a_F^{II})\cos{w}+b_F^{II}
\hspace{4 mm}
(c_F^{II}-a_F^{II})\sin{w}
  \end{array} \}^T
\end{eqnarray}
The solution to the set of equations (\ref{25}) gives the real parts of the numerical wavenumbers in the form:

\begin{eqnarray}\label{29}
\Re{(w_F')} = 
\frac{\sin{w}[f_1^I+f_2^I\cos{w}+f_3^I\cos{2w}]}
{g_1+g_2\cos{w}+g_3\cos{2w}}
\end{eqnarray}
and

\begin{eqnarray}\label{30}
\Re{((w_F'')^2)} = 
\frac{f_1^{II}+f_2^{II}\cos{w}+f_3^{II}\cos{2w}+f_4^{II}\cos{3w}}
{g_1+g_2\cos{w}+g_3\cos{2w}}
\end{eqnarray}
where $f_1^I$, $f_2^I$, $f_3^I$, $g_1$, $g_2$, $g_3$, $f_1^{II}$, $f_2^{II}$, $f_3^{II}$, $f_4^{II}$ are functions of the weights in equations (\ref{21}) and (\ref{22}). Comparing the equations (\ref{29}) and (\ref{30}) to the equations (\ref{19}) and (\ref{20}), the weights in equations (\ref{21}) and (\ref{22}) can be determined by solving the next polynomial equations:

\begin{eqnarray}\label{} \nonumber
\beta_F^Ic_F^{II}\theta_F^I + a_F^I + \beta_F^Ib_F^I - c_F^I - \theta_F^Ia_F^I\beta_F^{II} \\ \nonumber
- \theta_F^Ic_F^I\beta_F^{II} + c_F^{II}\theta_F^I\theta_F^{II} + 2\beta_F^Ia_F^I\theta_F^{II} - \theta_F^Ib_F^I\beta_F^{II}\theta_F^{II} + a_F^I\theta_F^{II2} \\
+ \beta_F^Ib_F^I\theta_F^{II2} - c_F^I\theta_F^{II2} - \beta_F^I\theta_F^Ia_F^{II} - \theta_F^I\theta_F^{II}a_F^{II} \\ \nonumber
+ \theta_F^Ib_F^{II} + \theta_F^{I2}\beta_F^{II}b_F^{II} - \beta_F^I\theta_F^I\theta_F^{II}b_F^{II} = -293/216; \nonumber
\end{eqnarray}

\begin{eqnarray}\label{}  \nonumber
c_F^{II}\theta_F^I - \theta_F^Ib_F^I\beta_F^{II} + a_F^I\theta_F^{II} - c_F^I\theta_F^{II} \\
- \theta_F^Ia_F^I\beta_F^{II}\theta_F^{II} + \theta_F^{I2}\beta_F^{II}a_F^{II} + \beta_F^I(a_F^I \\ \nonumber
+ b_F^I\theta_F^{II} + a_F^I\theta_F^{II2} - \theta_F^I\theta_F^{II}a_F^{II}) = -63/216; \nonumber
\end{eqnarray}

\begin{eqnarray}\label{} 
- 2a_F^I(\theta_F^I\beta_F^{II} - \beta_F^I\theta_F^{II}) = -1/216;
\end{eqnarray}

\begin{eqnarray}\label{} 
1 + \beta_F^{I2} + \theta_F^{I2}\beta_F^{II2} + 2\beta_F^I\theta_F^{II} - 2\beta_F^I\theta_F^I\beta_F^{II}\theta_F^{II} + \theta_F^{II2} + \beta_F^{I2}\theta_F^{II2} = 34/36;
\end{eqnarray}

\begin{eqnarray}\label{} 
2(\beta_F^I + \theta_F^{II})(1 - \theta_F^I\beta_F^{II} + \beta_F^I\theta_F^{II}) = 11/12;
\end{eqnarray}

\begin{eqnarray}\label{} 
-2\theta_F^I\beta_F^{II} + 2\beta_F^I\theta_F^{II} = 1/12;
\end{eqnarray}

\begin{eqnarray}\label{}  \nonumber
\beta_F^Ic_F^{II} - \beta_F^Ic_F^{II}\theta_F^I\beta_F^{II} - a_F^I\beta_F^{II} - \beta_F^Ib_F^I\beta_F^{II} \\ 
+ \theta_F^Ic_F^I\beta_F^{II2} + c_F^{II}\theta_F^{II} + \beta_F^{I2}c_F^{II}\theta_F^{II} - b_F^I\beta_F^{II}\theta_F^{II} - \beta_F^Ic_F^I\beta_F^{II}\theta_F^{II} + \beta_F^Ia_F^{II} \\ \nonumber
+ b_F^{II} + \beta_F^{I2}b_F^{II} + \beta_F^I\theta_F^{II}b_F^{II} = -1730/1296; \nonumber
\end{eqnarray}

\begin{eqnarray}\label{}  \nonumber
-b_F^I\beta_F^{II} + \theta_F^Ib_F^I\beta_F^{II2} - a_F^I\beta_F^{II}\theta_F^{II} - c_F^I\beta_F^{II}\theta_F^{II} \\ 
+ c_F^{II}(1 + \beta_F^{I2} - \theta_F^I\beta_F^{II} + 2\beta_F^I\theta_F^{II}) + a_F^{II} + \theta_F^{II}b_F^{II} - \beta_F^I(a_F^I\beta_F^{II} + c_F^I\beta_F^{II} \\ \nonumber
+ b_F^I\beta_F^{II}\theta_F^{II} - \theta_F^{II}a_F^{II} - 2b_F^{II} + \theta_F^I\beta_F^{II}b_F^{II}) + \beta_F^{I2}(a_F^{II} + \theta_F^{II}b_F^{II}) = 675/1296; \nonumber
\end{eqnarray}

\begin{eqnarray}\label{}  \nonumber
-c_F^I\beta_F^{II} + \theta_F^Ia_F^I\beta_F^{II2} + \theta_F^{II}a_F^{II} + \beta_F^{I2}\theta_F^{II}a_F^{II} \\
- \theta_F^I\beta_F^{II}b_F^{II} + \beta_F^I(c_F^{II} - a_F^I\beta_F^{II}\theta_F^{II} \\ \nonumber
+ a_F^{II} - \theta_F^I\beta_F^{II}a_F^{II} + \theta_F^{II}b_F^{II}) = 10870/1296; \nonumber
\end{eqnarray}

\begin{eqnarray}\label{} 
-\theta_F^I\beta_F^{II}a_F^{II} + \beta_F^I\theta_F^{II}a_F^{II} = 29/1296
\end{eqnarray}

The next open problem is to solve the above closed set of polynomial equations by an iterative method to obtain the weights of the scheme in (\ref{21}) and (\ref{22}).



\end{document}